\documentclass{article}
\usepackage[top=3cm,bottom=3cm,left=3cm,right=3cm]{geometry}
\usepackage[dvipsnames,svgnames,x11names]{xcolor}
\usepackage{amsmath,amssymb,mathtools}
\usepackage{graphicx}
\usepackage{float}
\usepackage{booktabs,tabularx,multirow,colortbl}
\usepackage{caption}
\usepackage{hyperref}
\usepackage{fancyhdr}
\usepackage{enumitem}

\fancypagestyle{firstpaper}{
  \fancyhf{}
  
  \lfoot{\small \copyright\  Hassan Jolany}
}
\usepackage{multicol}
\usepackage{setspace}
\usepackage{xurl}

\begin{document}

\thispagestyle{empty}


\vspace{3cm}


\vspace{2cm}



\vspace{2cm}


\vspace*{\fill}

\hfill

\vspace{1cm}


\thispagestyle{firstpaper}
\setcounter{page}{1}

\noindent{\Large\bfseries Generalized Sasaki-Einstein metric twisted with Weil-Petersson metric}

\vspace{0.5em}
\noindent{\large\bfseries Hassan Jolany}

\vspace{1.5cm}

\begin{center}
\begin{minipage}{0.9\textwidth}
\small
Existence of canonical metric on a Sasakian variety was a long standing conjecture and the major part of this conjecture is about varieties which do not have definite basic first Chern class(most of the Sasakian varieties do not have definite basic first Chern class). We give a program by using Fujino's Minimal Model program to find canonical metric on the canonical model of Sasakian varieties.
\end{minipage}
\end{center}

\vspace{1.5cm}

\section{Introduction}

Kähler geometry study even dimensional manifolds and Sasakian Geometry is the odd dimensional analogue of Kähler geometry and special case of contact structure. Sasaki geometry in dimension $2n + 1$ is related to Kähler geometry in both dimension $2n$ and $2n + 2$. Sasakian manifolds are foliated manifolds with a contact structure, which are the generalization of Kähler manifolds. For details see [1], [2]. The geometry of Sasakian manifolds has recently garnered a great deal of interest due to its important role in the anti-de Sitter/conformal field theory correspondence of theoretical physic. In particular, a conformal field theory is dual to $AdS_5 \times M^5$, where $AdS_5$ is anti-de Sitter space and $M$ is a 5-dimensional Sasaki-Einstein manifold. Sasaki-Einstein manifolds are of interest in the study of non-compact Calabi-Yau manifolds.

\textbf{Definition 1.1.} A Riemannian manifold $(M, g)$ is Sasakian if the metric cone $(C(M), \bar{g})$, $C(M) := \mathbb{R}^+ \times M$ and $\bar{g} = dr^2 + r^2g$, is Kähler, i.e. $\bar{g}$ admits a compatible almost complex structure $J$ so that $(C(M), \bar{g}, J)$ is a Kähler structure. We embed $M$ in $C(M)$ as the level set $\{r = 1\}$. \hfill $\square$

Sasakian manifold $(M, g)$ inherits a contact structure, and Reeb vector field and so we can introduce Sasakian manifolds with the following new datas.

\begin{enumerate}[label=\bfseries \arabic*., leftmargin=1.5cm]
\item A contact structure
\begin{equation*}
\eta = d^c \log r = Jd \log r
\end{equation*}
\item With Reeb vector field
\begin{equation*}
\xi = Jr\partial_r \in \Gamma(TM)
\end{equation*}
where $\eta(\xi) = 1$, $i_\xi d\eta = 0$
\item A strictly pseudoconvex CR-structure $(D, I)$, $D = \ker \eta$.
\item $I$ induces a transversely holomorphic structure on $\mathcal{F}_\xi$, the Reeb foliation, with transverse Kähler form $\omega^T = \frac{1}{2}d\eta$
\end{enumerate}

There is an orthogonal decomposition of the tangent bundle
\begin{equation*}
TM = D \oplus L_\xi
\end{equation*}
where $L_\xi$ is the trivial bundle generalized by $\xi$. We can then introduce the $(1, 1)$-tensor field $\Phi$ such that
\begin{equation*}
\Phi(\xi) = 0 \text{ and } \Phi(X) = JX \text{ on } X \in \Gamma(D)
\end{equation*}
It holds that $\Phi^2 = -\mathrm{id} + \xi \otimes \eta$ and $g(\Phi\cdot, \Phi\cdot) = g + \eta \otimes \eta$

Note that $\Phi$ is compatible with the 2-form $d\eta$,
\begin{equation*}
d\eta(\Phi X, \Phi Y) = d\eta(X, Y), \ X, Y \in \Gamma(TM)
\end{equation*}


and $g(X,Y) = \frac{1}{2} d\eta (X,\Phi Y), X,Y\in \Gamma (D)$. Hence $d\eta$ defines a symplectic form on $D$.

Since $\mathcal{L}_{r\partial_r}\omega = 2\omega$, we can write globally,
\begin{equation*}
\omega = \frac{1}{2} \sqrt{-1} \partial \bar{\partial} r^2
\end{equation*}
where $\omega$ is the Kähler form associated to $\bar{g}$ on $C(M)$. Hence $\frac{1}{2} r^2$ is a global Kähler potential on $C(M)$.

The cone $\mathcal{M} := C(M) \cup \{o\}$ of a Sasakian manifold $M$ is an affine algebraic variety with an algebraic action of some torus $T^r = (\mathbb{C}^*)^r$ and can be embedded as $\mathcal{M} \hookrightarrow \mathbb{C}^N$.

The canonical example of a Sasaki-Einstein manifold is the odd dimensional sphere $S^{2n - 1}$, equipped with its standard Einstein metric. In this case the Kähler cone is $\mathbb{C}^n \setminus \{0\}$, equipped with its flat metric.

Sasakian manifolds can be characterized into three categories based on the orbits of the Reeb field. If the orbits of the Reeb field are all closed, then they generate a locally free, isometric $U(1)$ action on $(M,g)$. If the $U(1)$ action is free, then $(M,g)$ is said to be regular and the quotient manifold $M / U(1)$ is Kähler. If the action is not free, then $(M,g)$ is said to be quasi-regular, and the quotient manifold is a Kähler orbifold. If the orbits are not closed, then the Sasakian manifold $(M,g)$ is said to be irregular. It was widely believed that irregular Sasaki-Einstein Manifolds did not exist as per a conjecture by Tian and Cheeger. However, in 2004 a landmark paper of Sparks, Martelli, Gauntlett and Waldram [7], constructed an infinite sequence of Sasaki-Einstein metrics on $S^2\times S^3$ which included irregular and quasiregular Sasaki-Einstein Manifolds.

We need to recall the notion of transverse cohomology.

\textbf{Definition 1.2.} A $p$-form $\omega$ on $(M, \xi)$ is called basic if $i_{\xi}\omega = 0$, and $\mathcal{L}_{\xi}\omega = 0$. The sheaf of basic functions is denoted by $C_B^\infty(M)$. On the Sasakian manifold $(M, \xi)$ we take $\Lambda_B^p$ be the sheaf of basic $p$-forms, and $\Omega_B^p = \Gamma(M, \Lambda_B^p)$ the global sections. The de Rham differential $d$ preserves basic forms, and hence restricts to a well defined operator $d_B: \Lambda_B^p \to \Lambda_B^{p+1}$. We thus get a complex
\begin{equation*}
0 \to C_B^\infty(M) \to \Omega_B^1 \xrightarrow{d_B} \dots \xrightarrow{d_B} \Omega_B^{2n} \xrightarrow{d_B} 0
\end{equation*}
The basic de Rham cohomology groups denoted by $H_B^p(M)$. \hfill $\square$

By transverse complex structure $\Phi$ we can decompose
\begin{equation*}
\Lambda_B^r \otimes \mathbb{C} = \bigoplus_{p+q=r} \Lambda_B^{p,q}
\end{equation*}
and hence we can write $d_B = \partial_B + \bar{\partial}_B$, and $d_B^c = \frac{1}{2}\sqrt{-1}(\bar{\partial}_B - \partial_B)$ where
\begin{equation*}
\partial_B: \Lambda_B^{p,q} \to \Lambda_B^{p+1,q} \text{ and } \bar{\partial}_B: \Lambda_B^{p,q} \to \Lambda_B^{p,q+1}
\end{equation*}
and it is clear that $d_B^2 = d_B d_B^c = 0$.

The transverse metric $g^T$ is related to the Sasaki metric $g$ by
\begin{equation*}
g = g^T + \eta \otimes \eta
\end{equation*}
From the transverse metric $g^T$, we can define the transverse Levi-Civita connection $\nabla^T$ on $D$ by
\begin{equation*}
\nabla_X^T Y = (\nabla_X Y)^p, \ X, Y \in \Gamma(D)
\end{equation*}
\begin{equation*}
\nabla_\xi^T Y = [\xi, Y]^p, \ Y \in \Gamma(D)
\end{equation*}
where $X^p$ denotes the projection of $X$ onto $D$. The transverse Levi-Civita connection is the unique torsion free connection.

\textbf{Definition 1.3.} The space of Sasaki metrics is defined as
\begin{equation*}
\mathcal{H} = \{\phi \in C_B^\infty(M) \mid \eta_\phi \wedge (d\eta_\phi)^n \neq 0, \eta_\phi = \eta + d_B^c \phi\}
\end{equation*}
\hfill $\square$


\textbf{Definition 1.4.} We can then define the transverse curvature operator by
\begin{equation*}
Rm^T(X, Y) = \nabla_X^T \nabla_Y^T - \nabla_Y^T \nabla_X^T - \nabla_{[X, Y]}^T
\end{equation*}
and by the same way we can define the transverse Ricci curvature and transverse scalar curvature. Take the form $\rho^T = Ric^T(\Phi \cdot, \cdot)$, is called the transverse Ricci form. One can see that
\begin{equation*}
\rho^T = -\sqrt{-1}\partial_B \bar{\partial}_B \log \det(g^T)
\end{equation*}
The $(1,1)$-basic Dolbeault cohomology class $[\rho^T]_B \in H_B^{1,1}(M)$ is called the basic first Chern class which is independent of the transverse metric $g^T$ and denoted by $c_1^B(M)$. \hfill $\square$

Now we recall the transverse Einstein metric which is an analogue of the Kähler-Einstein metric in Kähler setting.

\textbf{Definition 1.5.} A transverse metric $g^T$ is called a transverse Einstein metric if it satisfies in $Ric^T = cg^T$ for some constant $c$. The transverse metric $g^T$ is said to have transverse constant scalar curvature if
\begin{equation*}
Tr_{g^T} Ric^T = c
\end{equation*}
\hfill $\square$

\textbf{Theorem 1.6.} The metric $g$ is Sasaki-Einstein metric if and only if $g^T$ is Kähler-Einstein if and only if $\bar{g}$ is Ricci-flat metric. Moreover, If there exists a transverse Kähler Einstein metric with $Ric^T = \lambda \rho^T$ for some constant $\lambda$, then $c_1^B = \lambda [d\eta]_B \in H_B^{1,1}(M)$. This implies that $c_1^B$ definite depending on the sign of $\lambda$ and $c_1(D) = 0$. It immediately follows that for a Sasaki-Einstein manifold the restricted holonomy group of the cone $Hol^0(\bar{g}) \subset SU(n)$. \hfill $\square$

\section{Song-Tian program on Sasakian manifolds}

El Kacimi-Alaoui [3] showed that by assuming $c_1^B = \lambda [d\eta]_B$, there is a basic function $F$ such that
\begin{equation*}
\rho^T - \lambda d\eta = \sqrt{-1}\partial_B \bar{\partial}_B F
\end{equation*}
The transverse Kähler-Einstein equation can be reduced to the following transverse Monge-Ampere equation
\begin{equation*}
\frac{\det(g_{i\bar{j}}^T + \phi_{i\bar{j}})}{\det(g_{i\bar{j}}^T)} = e^{-\lambda \phi + F}, \ g_{i\bar{j}}^T + \phi_{i\bar{j}} > 0
\end{equation*}
where here $\phi$ is basic. Note that, it is not elliptic, but transversal elliptic.

\textbf{Sasakian Calabi Problem(Boyer-Galicki):} Give a manifold $M$ with Sasakian structure $(\xi, \eta, \Phi, g)$ and the basic first Chern class $c_1^B$ is positive, negative or null, can one deform it to another Sasakian structure $(\xi, \eta', \Phi', g')$ with an $\eta$-Einstein metric $g'$?

We have the same analogue of Kähler Ricci flow in Sasakian setting and is called Sasaki-Ricci flow [4].

\textbf{Theorem 2.1.} On a compact manifold with Sasakian structure $(M,\xi ,\eta ,\Phi ,g)$, $c_{1}^{B} = \lambda [d\eta ]_{B}$. There is a smooth family of Sasakian structures $(\xi ,\eta (t),\Phi (t),g(t))$ satisfying $(\xi ,\eta (0),\Phi (0),g(0)) = (\xi ,\eta ,\Phi ,g)$ and
\begin{equation*}
\frac{\partial}{\partial t} g^T(t) = -\left(Ric^T(g(t)) - \lambda g^T(t)\right)
\end{equation*}
and we can write it as transverse Monge-Ampere equation
\begin{equation*}
\frac{\partial}{\partial t} \varphi = \log \det(g_{i\bar{j}}^T + \varphi_{i\bar{j}}) - \log \det(g_{i\bar{j}}^T) + \lambda \varphi - F
\end{equation*}
When the basic first Chern class $c_1^B$ is negative or null, then the Sasaki Ricci flow converges to $\eta$-Einstein metric, [4]. When the basic first Chern class $c_1^B$ is positive we need the algebro-geometric notion of K-stability in Sasakian setting which recently solved by T. Collins and Gábor Székelyhidi [8]. When the Sasakian manifold is quasi-regular, then this is equivalent to the work of Ross-Thomas on K-stability for orbifolds. \hfill $\square$


But the fact is that the most of Sasaki manifolds do not have definite basic first Chern class and the question is, how can we find the canonical metric (generalized Sasaki-Einstein metric) on it.

The fact is that Sasakian manifold $M$ is in deal with affine variety $\mathcal{M}$ instead of projective variety and we must run MMP for affine variety to get canonical metric. We first need to introduce canonical model of $\mathcal{M}$ to get canonical metric.

Let $X_0$ be a projective variety with canonical line bundle $K \to X_0$ of Kodaira dimension
\begin{equation*}
\kappa(X_0) = \limsup \frac{\log \dim H^0(X_0, K^{\otimes \ell})}{\log \ell}
\end{equation*}
This can be shown to coincide with the maximal complex dimension of the image of $X_0$ under pluri-canonical maps to complex projective space, so that $\kappa(X_0) \in \{-\infty, 0, 1, \dots, m\}$. Also since in general we work on Singular Kähler variety we need to notion of numerical Kodaira dimension instead of Kodaira dimension.
\begin{equation*}
\kappa_{num}(X) = \sup_{k \geq 1} \left[ \limsup_{m \to \infty} \frac{\log \dim_{\mathbb{C}} H^0(X, mK_X + kL)}{\log m} \right]
\end{equation*}
where $L$ is an ample line bundle on $X$. Note that the definition of $\kappa_{num}(X)$ is independent of the choice of the ample line bundle $L$ on $X$. Siu formulated that the abundance conjecture is equivalent as
\begin{equation*}
\kappa_{kod}(X) = \kappa_{num}(X)
\end{equation*}
Numerical dimension is good thing.

It is worth to mention that if $f: X \to Y$ be an algebraic fibre space and $\kappa(X) \geq 0$, $\kappa(Y) = \dim Y$, (for example Iitaka fibration), then $\kappa(X) = \kappa(Y) + \kappa(F)$, where $F$ is a general fibre of $F$.

Let $V$ be an irreducible algebraic variety. By Nagata's theorem, we have a complete algebraic variety $\bar{V}$ which contains $V$ as a dense Zariski open subset. By Hironaka's theorem, we have a smooth projective variety $\bar{W}$ and a projective birational morphism $\mu : \bar{W} \to \bar{V}$ such that if $W = \mu^{-1}(V)$ then $\bar{D} = \bar{W} - W = \mu^{-1}(\bar{V} - V)$, is a simple normal crossing divisor on $\bar{W}$. The logarithmic Kodaira dimension $\bar{\kappa}(V)$ of $V$ is defined as
\begin{equation*}
\bar{\kappa}(V) = \kappa(\bar{W}, K_{\bar{W}} + \bar{D})
\end{equation*}
$\bar{\kappa}(V)$ is well-defined, that is, it is independent of the choice of the pair $(\bar{W},\bar{D})$. We have the following theorem from Yoshinori Gongyo and Osamu Fujino [6].

\textbf{Theorem 2.2.} Let $V$ be an affine variety. We can take a pair $(\bar{W}, \bar{D})$ such that
\begin{equation*}
\bar{\kappa}(V) = \kappa(\bar{W}, K_{\bar{W}} + \bar{D})
\end{equation*}
More precisely, by running a minimal model program with ample scaling, we have a finite sequence of flips and divisorial contraction
\begin{equation*}
(\bar{W}, \bar{D}) = (\bar{W}_0, \bar{D}_0) \dashrightarrow (\bar{W}_1, \bar{D}_1) \dashrightarrow \dots \dashrightarrow (\bar{W}_k, \bar{D}_k)
\end{equation*}
such that $(\bar{W}_k,\bar{D}_k)$ is a good minimal model or has a Mori fiber space structure. \hfill $\square$

Now we introduce the canonical model of affine variety $V$.

\textbf{Definition 2.3.} Let $V$ be an affine variety and let $(\bar{W},\bar{D})$ be a pair as in previous theorem. We define the canonical model of $V$ as
\begin{equation*}
V_{can} = \operatorname{Proj} \bigoplus_{m \geq 0} H^0(\bar{W}, \mathcal{O}_{\bar{W}}(m(K_{\bar{W}} + \bar{D})))
\end{equation*}
and the $g_{can}$ is the canonical metric on the canonical model of $V$. \hfill $\square$

Abundance conjecture tells us that if a minimal model exists, then the canonical line bundle $V_{min}$ induces a unique holomorphic map
\begin{equation*}
\pi : V_{min} \to V_{can}
\end{equation*}
where $V_{can}$ is the unique canonical model of $V_{min}$.

It is easy to see that the canonical ring $R(V) = \bigoplus_{m\geq 0}H^0 (\bar{W},\mathcal{O}_{\bar{W}}(m(K_{\bar{W}} + \bar{D})))$ is independent of the pair $(\bar{W},\bar{D})$ and is well-defined. Then $R(V)$ is a finitely generated $\mathbb{C}$-algebra. This is because we can choose $(\bar{W},\bar{D})$ such that it has a good minimal model or a Mori fiber space structure.

By extending Song-Tian program[10, 11] on Sasakian manifold $M$, it turns out that the normalized Sasaki-Ricci flow doing exactly same thing to replace the Sasakian cone $\mathcal{M}$ by its minimal model by using finitely many geometric surgeries and then deform minimal model to canonical model.


\noindent\textbf{Theorem 2.4.} Let $M$ be a Sasakian manifold and $\mathcal{M}$ with $\bar{\kappa}(\mathcal{M}) \geq 0$ be the affine Sasakian cone, since the canonical ring $R(\mathcal{M})$ is finitely generated C-algebra then from [9],
\begin{equation*}
\pi : (\overline {{\mathcal {M}}}, \bar {D}) \to \mathcal {M} _ {c a n} ^ {D}
\end{equation*}
and we have
\begin{equation*}
Ric(g_{can}) = -g_{can} + g_{WP}^{\bar{D}} + + \sum_{P} (b(1 - t_P^{\bar{D}})) [\pi^*(P)] + [B^{\bar{D}}]
\end{equation*}
where $B^{\bar{D}}$ is $\mathbb{Q}$-divisor on $X$ such that $\pi_{*}\mathcal{O}_{X}([iB^{\bar{D}}_{+}])=\mathcal{O}_{B}\left(\forall i>0\right)$. Here $s_{P}^{\bar{D}}:=b(1-t_{P}^{\bar{D}})$ where $t_{P}^{\bar{D}}$ is the log-canonical threshold of $\pi^{*}P$ with respect to $(X,\bar{D}-B^{\bar{D}}/b)$ over the generic point $\eta_{P}$ of $P$. i.e.,
\begin{equation*}
t_{P}^{\bar{D}} := \max \{t \in \mathbb{R} \mid \left(X, \bar{D} - B^{\bar{D}}/b + t\pi^*(P)\right) \text{ is sub log canonical over } \eta_P\} \hfill \square
\end{equation*}

\vspace{1em}
Sasakian manifolds are equipped with a contact CR-structure and a CR-holomorphic action of the corresponding Reeb field. An analogue of a holomorphic map is obviously a CR-holomorphic submerssion $X \rightarrow B$ . One would also require that the Reeb field action on $B$ preserve $X$.

\vspace{1em}
\noindent\textbf{Theorem 2.5.} Let $X$ be a Sasakian manifold and $\pi : X \to B$ is a CR-holomorphic submerssion to a compact Sasakian manifold $B$ with $c_1^B(K_B) < 0$ where the general fibers are Calabi-Yau manifolds, i.e., $c_1(K_{X_s}) = 0$, and central fiber is Calabi-Yau manifold, Then $X$ admits a unique smooth transverse Kähler-Einstein metric $\omega_B$ solving
\begin{equation*}
Ric^T(\omega_B) = -\omega_B^T + \pi^*\omega_{WP}
\end{equation*}
where $\omega_{WP}$ is the Weil-Petersson form on the moduli space of Calabi-Yau fibers. \hfill $\square$

\vspace{1.5em}
\noindent\textbf{Proof.} Set $\pi^*(\omega_B^T) = -\sqrt{-1}\partial \bar{\partial}\log \Omega$ where $\Omega$ is a relative volume form on $D$.

Since fibers are Calabi-Yau manifolds so $c_1(X_y) = 0$, hence there is a smooth function $F_y$ such that $Ric(\omega_y) = \sqrt{-1}\partial \bar{\partial} F_y$ and $\int_{X_y}(e^{F_y} - 1)\omega_y^{n - m} = 0$. The function $F_y$ vary smoothly in $y$. By Yau's theorem there is a unique Ricci-flat Kähler metric $\omega_{SRF,y}$ on $X_y$ cohomologous to $\omega_0$. So there is a smooth basic function $\rho_y$ on $\pi^{-1}(y) = X_y$ such that $\omega_0|_{X_y} + \sqrt{-1}\partial \bar{\partial}\rho_y = \omega_{SRF,y}$ is the unique Ricci-flat Kähler metric on $X_y$. If we normalize by $\int_{X_y}\rho_y(\omega_0)^n |_{X_y} = 0$ then $\rho_y$ varies in $y$ and defines a function $\rho$ on $X$ and we let
\begin{equation*}
\omega_{SRF}^T = \omega_0^T + \sqrt{-1}\partial \bar{\partial}\rho
\end{equation*}
which is called as transverse Semi-Ricci Flat metric. Such Semi-Flat Calabi-Yau metrics were first constructed by Greene-Shapere-Vafa-Yau on surfaces in Kähler setting. More precisely, a closed real $(1,1)$-form $\omega_{SRF}$ on open set $U \subset X \setminus S$, (where $S$ is proper analytic subvariety contains singular points of $X$) will be called semi-Ricci flat if its restriction to each fiber $X_y \cap U$ with $y \in f(U)$ be Ricci-flat. Notice that $\omega_{SRF}$ is semi-positive

We have
\begin{equation*}
\pi_*(\Omega) = \frac{\Omega}{(\omega_{SRF})^m} = \int_{X_s} \Omega
\end{equation*}
Now solve the complex Monge-Ampere equation on $B$
\begin{equation*}
\left((1 - e^{-t})\omega_B^T + e^{-t}\omega_{WP} + \sqrt{-1}\partial \bar{\partial}v^B\right)^m = e^{v^B}\pi_*\Omega
\end{equation*}
Take, $\omega = \omega(t) = \omega_B^T + \sqrt{-1}\partial \bar{\partial}v^B$, we have
\begin{align*}
Ric^T(\omega) &= \\
&= -\sqrt{-1}\partial_B \bar{\partial}_B \log(\omega^T)^m \\
&= -\sqrt{-1}\partial_B \bar{\partial}_B \log\pi_*\Omega - \sqrt{-1}\partial_B \bar{\partial}_B v^B
\end{align*}
and

\begin{align*}
\sqrt{-1} \partial_B \bar{\partial}_B \log \pi_* \Omega + \sqrt{-1} \partial_B \bar{\partial}_B v^B &= \\
&= \sqrt{-1} \partial_B \bar{\partial}_B \log \pi_* \Omega + \omega^T - \omega_B^T
\end{align*}
Hence, by using $\pi_*(\Omega) = \frac{\Omega}{(\omega_{SRF})^m}$, we get
\begin{align*}
\sqrt{-1} \partial_B \bar{\partial}_B \log \pi_* \Omega + \sqrt{-1} \partial_B \bar{\partial}_B v^B &= \\
&= \omega^T - \pi^* \omega_{WP}
\end{align*}
So,
\begin{equation*}
Ric^T(\omega) = -\omega^T + \pi^* \omega_{WP}
\end{equation*}
where the Weil-Petersson metric is $\omega_{WP} = Ric\pi_*(K_{X/B})$

We consider the following normalized transverse-Ricci flow.
\begin{equation*}
\frac{\partial \omega_t^T}{\partial t} = -Ric^T(\omega_t) - \omega_t^T
\end{equation*}
with any transverse Kähler metric $\omega_0^T$ as the initial metric. We have $\pi^*[\omega_B^T] = -c_1^B(X)$ for some transverse Kähler class $[\omega_B^T]$ over $B$. Now take the reference metric
\begin{equation*}
\hat{\omega}_t = (1 - e^{-t})\pi^*(\omega_B^T) + e^{-t}\omega_0^T
\end{equation*}
Now $[\hat{\omega}_t^T] = [\omega_t^T]$ and from $\partial_B \bar{\partial}_B$-lemma we can write $\omega_t^T = \hat{\omega}_t^T + \sqrt{-1} \partial_B \bar{\partial}_B \varphi_t$ for some basic function $\varphi_t$. Let $\Omega$ be a volume form on $X$ such that
\begin{equation*}
\pi^*(\omega_B^T) = \sqrt{-1} \partial_B \bar{\partial}_B \Omega
\end{equation*}
The transverse Kähler Ricci flow is equivalent with the following Monge-Ampere equation
\begin{equation*}
\frac{\partial \varphi_t}{\partial t} = \log \frac{(\hat{\omega}_t^T + \sqrt{-1} \partial_B \bar{\partial}_B \varphi_t)^n}{e^{-rt} \Omega_{X/B}} - \varphi_t
\end{equation*}
From maximum principle we have $|\varphi_t| < C$. Also by simple computation we have
\begin{equation*}
\left( \frac{\partial}{\partial t} - \Delta_B \right) \left(-\varphi_t - (1 - e^{-t}) \frac{\partial \varphi_t}{\partial t}\right) = \Delta_B \left(-\varphi_t - (1 - e^{-t}) \frac{\partial \varphi_t}{\partial t}\right) - \text{Tr}_{\omega_t^T} \omega_0^T + re^{-t} + n - r
\end{equation*}
Then by applying Maximum Principle and the bounds for $\varphi_t$, we have an upper bound for $\frac{\partial \varphi_t}{\partial t} < C$. Now from the stimate of $|\varphi_t| < C$, $\frac{\partial \varphi_t}{\partial t} < C$, and using our Monge-Ampere equation, we have
\begin{equation*}
\left( \frac{\partial}{\partial t} - \Delta_B \right) \left(2\varphi_t + \frac{\partial \varphi_t}{\partial t}\right) > \frac{\partial \varphi_t}{\partial t} - C + Ce^{-\frac{\partial \varphi_t}{\partial t}}
\end{equation*}
and again by applying Maximum Principle, we have a lower bound for $\frac{\partial \varphi_t}{\partial t} > C'$. By using these estimates in our Monge-Ampere equation we have,
\begin{equation*}
C^{-1} e^{-rt} \Omega \le (\omega_t^T)^n \le Ce^{-rt} \Omega
\end{equation*}
where $\dim X_z = \pi^{-1}(z) = r$. Now, take the following Monge-Ampere equation
\begin{equation*}
(\omega_B^T + \sqrt{-1} \partial_B \bar{\partial}_B u)^{n-r} = Fe^u(\omega_B^T)^{n-r}
\end{equation*}
where
\begin{equation*}
F = \frac{\Omega}{\binom{n}{r} (\omega_B^T)^{n-r} \wedge (\omega_{SRF})^r}
\end{equation*}
The solution $\varphi_t$ for our Monge-Ampere equation converges uniformly to $u$ as $t \to \infty$. Now set $v_t = \varphi_t - u - e^{-t} \rho^B$. Then by taking time derivative of $v_t$ and rewriting it based on our Monge-Ampere equation


and using maximum principle we see that \( |\varphi_t - u| < Ce^{-t/2} \). Now by the same method of Song-Tian [1],[2], \( \hat{\omega}_t^T \leq \omega_t^T \leq C^{-1}\hat{\omega}_t^T \). Since \( \omega_t^T, \omega_0^T \) and \( \hat{\omega}_t^T \) are uniformly equivalent, we have

\begin{equation*}
| (\Delta_{B})_{\omega_{0}} \varphi_{t} | \leq C
\end{equation*}

and we get the \(C^{1,\alpha}\) convergence. \hfill \(\blacksquare\)

\vspace{1em}

Now we formulate G.Tian's conjecture in Sasakian setting via relative Sasaki Ricci flow. A compact Sasakian manifold \( X \) is said to admit a Fano fibration if there exists a CR-holomorphic immersion \( f: X \to Y \) with connected fibers and \( 0 \leq \dim Y < \dim X \) and such that \( -K_X \) is \( f \)-ample, i.e., with positive basic first Chern class of fibers \( c_1^B(X_y) > 0 \) and the generic fiber is of \( \dim X - \dim Y \). If \( Y \) be a point then we say \( X \) is Fano Sasakian manifold with positive basic first Chern class \( c_1^B(X) > 0 \).

\vspace{1em}

\textbf{Conjecture 1.} Let \( X^n \) be a compact Sasakian manifold. Then there exists a Sasaki metric \( \omega_0^T \) such that the relative Sasaki-Ricci flow

\begin{equation*}
\left\{ \begin{array}{c} \frac {\partial \omega^ {T} (t)}{\partial t} = - R i c _ {X / Y} ^ {T} (\omega (t)) \\ \omega^ {T} (0) = \omega_ {0} ^ {T} \end{array} \right.
\end{equation*}

has finite time collapsing if and only if \( X \) admits a Fano fibration \( f: X \to Y \). In this case, we can write

\begin{equation*}
[ \omega_ {0} ^ {T} ] = T (c _ {1} ^ {B} (X) + f ^ {*} (\omega_ {Y} ^ {T}))
\end{equation*}

for some Sasaki metric  \( \omega_{Y} \)  on Y, where T is the maximal existence time of the flow. \hfill \(\square\)

\vspace{1em}

\textbf{Conjecture 2.} Let \((X^n,\omega_0)\) be a compact Sasaki manifold, let \(\omega (t)\) be the solution of the Sasaki-Ricci flow, defined on the maximal time interval \([0,T)\) with \(T < \infty\).

\begin{equation*}
\left\{ \begin{array}{c} \frac {\partial \omega^ {T} (t)}{\partial t} = - R i c ^ {T} (\omega (t)) \\ \omega^ {T} (0) = \omega_ {0} ^ {T} \end{array} \right.
\end{equation*}

Then as \(t\to T\) we have

\begin{equation*}
\operatorname{diam} (X, \omega^ {T} (t)) \to 0
\end{equation*}

if and only if \([\omega_0^T] = \lambda c_1^B (X)\), for some \(\lambda > 0\).

\vspace{1em}

\textbf{Lemma 2.6.} For every Fano fibration for Sasakian manifold \( X \) there always exists a solution of the relative Sasaki-Ricci flow which collapses in finite time.

\vspace{1em}

\textbf{Proof.} Let we have Fano fibration. From the relative Sasaki-Ricci flow starting at  \( \omega_{0}^{T} \)

\begin{equation*}
\frac {\partial \omega^ {T} (t)}{\partial t} = - R i c _ {X / Y} ^ {T} (\omega (t))
\end{equation*}

from the definition of relative Ricci form \( Ric_{X/Y}^{T} \) for some Sasaki form \( \omega_{Y}^{T} \) in \( Y \), we get

\begin{equation*}
\frac {\partial [ \omega^ {T} (t) ]}{\partial t} = - c _ {1} ^ {B} (X) + f ^ {*} (\omega_ {Y} ^ {T})
\end{equation*}

hence

\begin{equation*}
[ \omega^ {T} (t) ] = [ \omega_ {0} ^ {T} ] + t (c _ {1} ^ {B} (K _ {X}) + f ^ {*} (\omega_ {Y} ^ {T}))
\end{equation*}

and since \( f \) is Fano fibration on Sasakian manifolds, we can write \( [\omega_0^T] = f^*(\omega_Y^T) + c_1^B(X) \) so

\begin{equation*}
[ \omega^ {T} (t) ] = (1 - t) (f ^ {*} (\omega_ {Y} ^ {T}) + c _ {1} ^ {B} (X))
\end{equation*}

showing that \([\omega ^T (t)]\) shrinks homothetically and would become degenerate at \(t = 1\). Moreover the total volume goes to zero when \(t\to 1\).


\section*{References}

\begin{enumerate}[label={[\arabic*]}, leftmargin=*, itemsep=0.5ex]
\item Boyer, C.P.; Galicki, K. On Sasakian-Einstein geometry. Internat. J. Math. 11(2000)873--909

\item C. Boyer, K. Galicki, and S.R.Simanca, Canonical Sasakian metrics, Commun. Math. Phys. 279(2008), 705--733

\item El Kacimi-Alaoui, A. Operateurs transversalement elliptiques sur un feuilletage riemannien et applications, Compositio Math. 79 (1990), 57--106

\item Smoczyk, Knut ;Guofang Wang, Zhang, Yongbing, The Sasaki-Ricci flow. Internat. J. Math. 21 (2010), no. 7, 951--969.

\item Tristan C. Collins, G{\'a}bor Sz{\'e}kelyhidi, Sasaki-Einstein metrics and K-stability, arXiv:1512.07213

\item Osamu Fujino, On subadditivity of the logarithmic Kodaira dimension, arXiv:1406.2759

\item Jerome P. Gauntlett, Dario Martelli, James Sparks, and Daniel Waldram, Sasaki-Einstein metrics on  \( S^{2} \times S^{3} \) , Adv. Theor. Math. Phys. 8 (2004), no. 4, 711--734. MR 2141499 (2006m:53067)

\item Tristan C. Collins, G{\'a}bor Sz{\'e}kelyhidi, Sasaki-Einstein metrics and K-stability, arXiv:1512.07213

\item Hassan Jolany, Canonical metric on canonical model of a mildly singular pair \((X;D)\)

\item Jian Song; Gang Tian, The K{\"a}hler-Ricci flow on surfaces of positive Kodaira dimension, Inventiones mathematicae. 170 (2007), no. 3, 609--653

\item Jian Song; Gang Tian, Canonical measures and K{\"a}hler-Ricci flow, J. Amer. Math. Soc. 25 (2012), no. 2, 303--353
\end{enumerate}

\end{document}